\documentclass[10pt,a4paper]{article}
\usepackage{amsfonts}
\usepackage{amssymb}
\usepackage{mathrsfs}
\usepackage{amsthm}
\usepackage{setspace}
\usepackage{color}
\usepackage{stmaryrd}

\usepackage{fancyhdr}
\newtheorem{theo}{Theorem}[section]
\newtheorem{prop}[theo]{Proposition}

\theoremstyle{definition}

\newcommand{\be}{\begin{eqnarray*}}
\newcommand{\ee}{\end{eqnarray*}}
\newcommand{\beqa}{\begin{eqnarray}}
\newcommand{\eeqa}{\end{eqnarray}}
\newcommand{\ba}{\begin{array}}
\newcommand{\ea}{\end{array}}

\newcommand{\mc}{\mathcal}
\newcommand{\mf}{\mathfrak}

\newcommand{\mbb}{\mathbb}

\begin{document}

\title{Addendum to ``Ricci-flat holonomy: a Classification'': the case of Spin(10)}
\author{Stuart Armstrong}
\date{April 2008}
\maketitle

\begin{abstract}
This note fills a hole in the author's previous paper ``Ricci-Flat Holonomy: a Classification'', by dealing with irreducible holonomy algebras that are subalgebras or real forms of $\mbb{C} \oplus \mf{spin}(10,\mbb{C})$. These all turn out to be of Ricci-type.
\end{abstract}

In my previous paper, \cite{meric}, I classified all possible irreducible holonomy algebras for torsion-free affine connections into three categories. I stated which holonomy algebras implied that their corresponding connections must be Ricci-flat, and gave a full list of those which were of Ricci-type (where the Ricci tensor encodes the full curvature tensor). But the list was incomplete...

One algebra family, namely $\mf{spin}(10)$ and the related algebras, was missing from the list. This short note will correct that omission. It will demonstrate that if $\nabla$ is a torsion-free affine connection with holonomy contained in $\mbb{C} \oplus \mf{spin}(10,\mbb{C})$ under one of the two standard $16$-dimensional representation, then $\nabla$ must be of Ricci type. This implies the same result for all subalgebras and real forms of $\mbb{C} \oplus \mf{spin}(10,\mbb{C})$, thus fixing the hole in the previous paper.

It will be usefull to recall some notation and notions from \cite{meric}. Given a Lie algebra $\mf{g}$ with a representation $V$, the formal curvature module $K(\mf{g})$ is the space of elements of $\wedge^2 V^* \otimes \mf{g}$ that obey the algebra\"ic Bianchi identity under the inclusion $\mf{g} \subset V^* \otimes V$.

Let $M$ be a manifold, and $\mc{G}$ a principal bundle over $M$ with structure group $G$ such that $TM = \mc{G} \times_G V$. Here $G$ is any Lie group with Lie algebra $\mf{g}$ and same representation $V$. Then if $\nabla$ is any torsion free connection with whose principal frame bundle reduces to a subbundle of $\mc{G}$, then $\nabla$'s curvature must take values in the bundle
\be
\mc{G} \times_G K(\mf{g}).
\ee
Consequently if $K(\mf{g})$ is of Ricci-type (in other words, every element of $K(\mf{g})$ is determined by its Ricci-like trace), then $\nabla$ is of Ricci-type.

From now on, let $\mf{g} = \mbb{C} \oplus \mf{spin}(10,\mbb{C})$ and $V$ be one of the two standard $16$-dimensional spin representations. The module $K(\mf{g})$ splits (see \cite{meric} and \cite{holclass}) as
\be
K(\mf{g}) = \partial (V^* \otimes \mf{g}^{(1)} ) \oplus H^{1,2}(\mf{g}).
\ee
Here, $\mf{g}^{(1)} = (V^* \otimes \mf{g}) \cap (\odot^2 (V^*)) \otimes V$, $\partial$ is the operator that acts on $\wedge^{k-1} V^* \otimes \mf{g}^{(k)}$ by anti-symmetrising the first $k$ $V^*$ terms, and $H^{1,2}$ is the (Spencer) cohomology component at $\wedge^2 V^* \otimes \mf{g}$ given by $\partial$.

Spencer co-homology is notoriously hard to calculate, however in this case we can work with Lie algebra comhomology for parabolic geometries \cite{TCPG}. In details, let $\mf{e}_6$ be the Lie algebra of the exceptional Lie group $E_6^{\mbb{C}}$. This algebra has a one-grading, i.e. it splits as
\be
\mf{e}_6 &=& \mbb{C}^{16} \oplus \big( \mbb{C} \oplus \mf{spin}(10,\mbb{C}) \big) \oplus \mbb{C}^{16*} \\
&=& V \oplus \mf{g} \oplus V^*.
\ee

There is a natural operator $\partial_{Lie}: \wedge^k V^* \otimes \mf{g} \to \wedge^{k+1} V^* \otimes \mf{g}$, given by:
\be
(\partial_{Lie} \psi)(v_0,\ldots v_{k}) &=& \sum_{i=0}^k (-1)^i [v_i, \psi(v_0, \ldots, \hat{v_i}, \ldots, v_k)],
\ee
where hats designate omission. This $\partial_{Lie}$ squares to zero and hence generates $H^k(V,\mf{g})$, the $k$-th Lie algebra co-homology of $(V,\mf{g})$.

\begin{prop}
On the module $\wedge^2 V^* \otimes \mf{g}$, $\partial = 3 \partial_{Lie}$. Therefore $K(\mf{g}) = H^{2}(V,\mf{g}) \cap (\wedge^2 V^* \otimes \mf{g}) \oplus \partial_{Lie} (V^* \otimes V)$.
\end{prop}
\begin{proof}
In $\mf{e}_6$, the Lie algebra bracket between $\mf{g}$ and $V$ is given by minus the standard action of $\mf{g}$ on $V$. Consequently, for $\psi$ an element of $\wedge^2 V^* \otimes \mf{g}$:
\be
(\partial_{Lie} \psi)(v_0,v_1,v_2) &=& \psi(v_1,v_2) \cdot v_0 - \psi(v_0,v_2) \cdot v_1 + \psi(v_0,v_1) \cdot v_2 \\
&=& 3 (\partial \psi)(v_0,v_1,v_2).
\ee
Now $\partial_{Lie}$ must respect homogeneity (see \cite{TCPG}), meaning that it must map $V^* \otimes V^*$ to $\wedge^2 V^* \otimes \mf{g}$, $V^* \otimes \mf{g}$ to $\wedge^2 V^* \otimes V$ and $V^* \otimes V$ to zero. Then since $\mf{g}^{(1)} = 0$,
\be
K(\mf{g}) &=& \ ker \ \partial = \ (ker \ \partial_{Lie})|_{\wedge^2 V^* \otimes \mf{g}} \\
&=& \big( H^{2}(V,\mf{g}) \cap (\wedge^2 V^* \otimes \mf{g})\big) \oplus \big(\partial_{Lie} (V^* \otimes V^*)\big).
\ee
\end{proof}

These equalities would be of little interest unless we could calculate the Lie algebra cohomologies. However, Kostant's version of the Bott-Borel-Weil theorem \cite{Kostant} allows one to do just that, giving the result:
\begin{theo}
$K(\mf{g}) \cap (\wedge^2 V^* \otimes \mf{g}) = 0$ and consequently $K(\mbb{C} \oplus \mf{spin}(10,\mbb{C}))$ is of Ricci-type.
\end{theo}
\begin{proof}
By Konstant's methods, we see that $H^2(V,\mf{g})$ is contained in $\wedge^2 V^* \otimes V$, meaning that $H^{2}(V,\mf{g}) \cap (\wedge^2 V^* \otimes \mf{g}) = 0$, and hence that $K(\mf{g}) = \partial_{Lie} (V^* \otimes V^*)$. However (see \cite{TCPG}), if $t_R$ is the Ricci-trace operator, then the compined operator $t_R \circ \partial_{Lie}$ is an automorphism of $V^* \otimes V^*$, meaning that any element of $K(\mf{g})$ is defined entirely by its Ricci trace.
\end{proof}

\bibliographystyle{alpha}
\bibliography{ref}

\end{document}